\documentclass{article}
\usepackage{amssymb,latexsym,amsmath,amsthm,mathrsfs}
\usepackage{graphicx}
\newcommand{\comment}[1]{}

\begin{document}
\title{A conjecture on the forms of the roots of equations\footnote{Presented to the St. Petersburg Academy on
November 2, 1733.
Originally published as
{\em De formis radicum aequationum cuiusque ordinis coniectatio},
Commentarii academiae scientiarum Petropolitanae \textbf{6} (1738),
216--231.
E30 in the Enestr{\"o}m index.
Translated from the Latin by Jordan Bell,
Department of Mathematics, University of Toronto, Toronto, Canada.
Email: jordan.bell@gmail.com}}
\author{Leonhard Euler}
\date{}
\maketitle

\S 1. It seems extremely remarkable that, although the roots
of cubic and biquadratic equations
were found
by means of the very first elements of analysis, still in these times,
when analysis
has grown by
leaps and bounds, the way still lies hidden for extracting the
roots
of higher equations, 
especially since
this matter was immediately investigated with the greatest diligence
by the most eminent savants.
From this diligence, even though the investigation has uncovered so little,
special methods for treating certain equations have been discovered.
For this reason I suppose that there will be no one who would censure
this undertaking of mine, in which I show what forms the roots of equations might
have and by what means they might perhaps be found,
even though I have not accomplished any more than that.
This may perhaps be helpful to lead others to finally solving this problem.

\S 2. Since an equation of any power includes in itself all lower powers,\footnote{Translator: Possibly Euler means that $\sum_{k=0}^m a_k x^k=0$ is a special
case of $\sum_{k=0}^n a_k x^k=0$ when $m<n$, taking $a_k=0$ for $m<k \leq n$.}
it is readily apparent that a method for extracting a root from such an equation
should be such that it involves methods for equations of all lower orders.
Hence the discovery of a root of an equation of the sixth degree cannot
be done unless the same has already been worked out for
equations of the fifth, fourth and third degree. Thus we 
see that the method of Bombelli for extracting the roots of
biquadratic equations leads to the resolution of the cubic equation;
and a root of the cubic equation cannot be defined without the resolution
of the quadratic equation.

\S 3. I will consider the resolution of the cubic equation
depending on the quadratic
in the following
way. Let the cubic equation
be
\[
x^3=ax+b,
\]
in which the second term is absent; I say that a root $x$ of it will be
\[
=\sqrt[3]{}A+\sqrt[3]{}B,
\]
with $A$ and $B$ being two roots of some quadratic equation
\[
z^2=\alpha z-\beta.
\]
Then from the nature of equations it will be
\[
A+B=\alpha, \quad AB=\beta.
\]
But for defining $\alpha$ and $\beta$ from $a$ and $b$, I take the equation
\[
x=\sqrt[3]{}A+\sqrt[3]{}B,
\]
which multiplied as a cube gives
\[
x^3=A+B+3\sqrt[3]{}AB(\sqrt[3]{}A+\sqrt[3]{}B)=3x\sqrt[3]{}AB+A+B.
\]
Comparing this with the given equation $x^3=ax+b$ will give
\[
a=3\sqrt[3]{}AB=3\sqrt[3]{}\beta \quad \textrm{and} \quad b=A+B=\alpha.
\]
It will thus become
\[
\alpha=b \quad \textrm{and} \quad \beta=\frac{a^3}{27};
\]
whence the quadratic equation serving for the resolution
of the equation $x^3=ax+b$ in the explained way will be
\[
z^2=bz-\frac{a^3}{27}.
\]
For with the roots $A$ and $B$ of this known, it will be
\[
x=\sqrt[3]{}A+\sqrt[3]{}B.
\]

\S 4. But since the cube root of any quantity has a triple value,
the formula $x=\sqrt[3]{}A+\sqrt[3]{}B$ will furthermore include all the roots of the
given
equation.
For let $\mu$ and $\nu$ be cube roots from unity besides unity;
it will further be
\[
x=\mu \sqrt[3]{}A+\nu \sqrt[3]{}B,
\]
if here it were $\mu \nu=1$. Consequently $\mu$ and $\nu$ must be 
\[
\frac{-1+\surd -3}{2} \quad \textrm{and} \quad \frac{-1-\surd -3}{2}
\]
or inversely. Therefore besides the root
\[
x=\sqrt[3]{}A+\sqrt[3]{}B
\]
the proposed equation is satisfied by these two other roots
\[
x=\frac{-1+\surd -3}{2}\sqrt[3]{}A+\frac{-1-\surd -3}{2}\sqrt[3]{}B
\]
and
\[
x=\frac{-1-\surd -3}{2}\sqrt[3]{}A+\frac{-1+\surd -3}{2}\sqrt[3]{}B.
\]
Of course, the roots of a cubic equation in which the second term is not
absent can also be determined by this rule.

\S 5. There are many different ways in which biquadratic equations are often
reduced to cubic equations which are however of no use to me.
 But I will use a particular method, like how before cubics were reduced
to quadratics, so that from this to a certain extent it can be concluded
how equations of higher degrees should be handled.
Thus if this equation is given
\[
x^4=ax^2+bx+c,
\]
in which again the second term is missing, I say it will be
\[
x=\surd A +\surd B+\surd C,
\]
where $A,B$ and $C$ are three roots of some cubic equation
\[
z^3=\alpha z^2-\beta z+\gamma.
\]
From this it will be
\[
\alpha=A+B+C, \quad \beta=AB+AC+BC \quad \textrm{and} \quad \gamma=ABC.
\]
To determine $\alpha,\beta$ and $\gamma$, the equation
$x=\surd A+\surd B+\surd C$ is freed from irrationality in this way.
The square is taken; it will be
\[
x^2=A+B+C+2\surd AB+2\surd AC+2\surd BC
\]
and then
\[
x^2-\alpha=2\surd AB+2\surd AC+2\surd BC.
\]
By taking the square again it will be
\[
x^4-2\alpha x^2+\alpha^2=4AB+4AC+4BC+8\surd ABC(\surd A+\surd B+\surd C)=
4\beta+8x\surd \gamma
\]
or
\[
x^4=2\alpha x^2+8x\surd \gamma+4\beta-\alpha^2.
\]
Comparing this equation with the given $x^4=ax^2+bx+c$ will give
\[
2\alpha=a, \quad 8\surd \gamma=b \quad \textrm{and} \quad 4\beta-\alpha^2=c,
\]
which yields
\[
\alpha=\frac{a}{2}, \quad \gamma=\frac{b^2}{64}, \quad \beta=\frac{c}{4}+
\frac{a^2}{16}.
\]
Therefore the cubic equation serving for the resolution of the biquadratic
equation is
\[
z^3=\frac{a}{2}z^2-\frac{4c+a^2}{16}z+\frac{b^2}{64}.
\]
For if the roots of this are $A,B$ and $C$, it will be
\[
x=\surd A+\surd B+\surd C.
\]
And the remaining three roots of the given equation will be
\[
\surd A-\surd B-\surd C,\quad \surd B-\surd A-\surd C \quad
\textrm{and} \quad \surd C-\surd A-\surd B.
\]

\S 6. Let us put $z=\surd t$; it will be
\[
\Big(t+\frac{4c+a^2}{16}\Big)\surd t=\frac{at}{2}+\frac{b^2}{64}
\]
and by taking the square we will have
\[
t^3+\frac{4c+a^2}{8}t^2+\frac{(4c+a^2)^2}{256}t=\frac{a^2t^2}{4}+
\frac{ab^2t}{64}+\frac{b^4}{4096}
\]
or
\[
t^3=\Big(\frac{a^2}{8}-\frac{c}{2}\Big)t^2+\Big(\frac{ab^2}{64}-\frac{cc}{16}-\frac{a^2c}{32}
-\frac{a^4}{256}\Big)t+\frac{b^4}{4096}.
\]
This equation therefore has the property that its roots are the squares
of the roots $A,B$ and $C$ of the prior equation. Whence if we put as the 
roots of this equation $E,F,G$, it will be
\[
x=\sqrt[4]{}E+\sqrt[4]{}F+\sqrt[4]{}G.
\]
And so a cubic equation is given, of which the biquadratic roots of its
roots, when summed, constitute a root of the given biquadratic
equation. And this method for finding the roots of biquadratic equations,
even if it is more laborious than the previous,
has a greater affinity with the resolution of cubic equations, since
a root of the same power as the proposed equation itself is extracted
from the roots of an inferior equation.

\S 7. Likewise by a similar rule the quadratic equation
\[
x^2=a,
\]
in which the second term is absent, can be resolved
by means of an equation of one dimension
\[
z=a.
\]
For since the root of this is $a$, and  hence the root
of the given equation
\[
x=\surd a \quad \textrm{or} \quad x=-\surd a.
\]
Indeed, the equation with inferior order, by means of which the 
superior equation missing a second term is resolved, I shall
call
the {\em resolvent equation}.
Thus the resolvent equation of the quadratic equation
\[
x^2=a
\]
will be
\[
z=a;
\]
the resolvent equation of the cubic equation
\[
x^3=ax+b
\]
will be
\[
z^2=bz-\frac{a^3}{27}
\]
and the resolvent equation of the biquadratic equation
\[
x^4=ax^3+bx+c
\]
will be
\[
z^3=\Big(\frac{a^2}{8}-\frac{c}{2}\Big)z^2-\Big(\frac{a^4}{256}+\frac{a^2c}{32}
+
\frac{c^2}{16}-\frac{ab^2}{64}\Big)z+\frac{b^4}{4096}.
\]
Namely, for the quadratic equation, if the root of the resolvent equation
is $A$, it will be
\[
x=\surd A;
\]
indeed for the cubic equation, if the roots of the resolvent are $A$ and $B$,
it will be
\[
x=\sqrt[3]{}A+\sqrt[3]{}B
\]
and too for the biquadratic equation, with the roots of the resolvent equation
being $A,B$ and $C$, it will be
\[
x=\sqrt[4]{}A+\sqrt[4]{}B+\sqrt[4]{}C.
\]

\S 8. From these cases, even if only three, it seems not to be without
sufficient reason for me to conclude that resolvent equations can be given
in this way.
Thus given the equation
\[
x^5=ax^3+bx^2+cx+d
\]
I infer that an equation of the fourth order is given
\[
z^4=\alpha z^3-\beta z^2+\gamma z-\delta,
\]
and that, if its roots are $A,B,C$ and $D$, then it will be
\[
x=\sqrt[5]{A}+\sqrt[5]{B}+\sqrt[5]{C}+\sqrt[5]{D}.
\]
And in general, the resolvent equation of the equation
\[
x^n=ax^{n-2}+bx^{n-3}+cx^{n-4}+\textrm{etc.},
\]
just as I suspect, will be
of the form
\[
z^{n-1}=\alpha z^{n-2}-\beta z^{n-3}+\gamma z^{n-4}-\textrm{etc.},
\]
whose roots are known and number $n-1$, say $A,B,C,D$ etc., and it will be
\[
x=\sqrt[n]{A}+\sqrt[n]{B}+\sqrt[n]{C}+\sqrt[n]{D}+\textrm{etc.}
\]
The truth of this conjecture could be consented to if resolvent equations
could be determined, for the roots of this equation can be quickly
assigned; namely by this progression
equations of lower order are continually led to, until
the true root of the given equation is disclosed.

\S 9. Although if the given equation has more than four dimensions
I am so far not able to define a resolvent equation,
however 
not insignificant evidence is at hand in which my conjecture is confirmed.
For if the given equation is such that
in the resolvent equation all the terms beside the first three vanish,
then this resolvent equation will always be able to be exhibited
and thus the roots of the given equation will be able to be defined.
In fact, the equations which admit resolution in this way are
exactly those which
the insightful Abraham de Moivre dealt with in
the Philosophical Transactions, no. 309. For 
if the resolvent equation were
\[
z^{n-1}=\alpha z^{n-2}-\beta z^{n-3}
\]
or
\[
z^2=\alpha z-\beta
\]
it will be possible to extract the equation that is to be resolved from this.
Let the roots of this equation be $A$ and $B$; for all the other
roots will vanish; a root of the equation that is to be resolved will be
\[
x=\sqrt[n]{}A+\sqrt[n]{}B.
\]
Indeed it is
\[
\alpha=A+B \quad \textrm{and} \quad \beta=AB
\]
from the nature of equations.
Thus it will then be
\[
\sqrt[n]{}A^2+\sqrt[n]{}B^2=x^2-2\sqrt[n]{}\beta
\]
and in turn
\begin{eqnarray*}
\sqrt[n]{}A^3+\sqrt[n]{}B^3&=&x^3-3x\sqrt[n]{}\beta,\\
\sqrt[n]{}A^4+\sqrt[n]{}B^4&=&x^4-4x^2\sqrt[n]{}\beta+2\sqrt[n]{}\beta^2,\\
\sqrt[n]{}A^5+\sqrt[n]{}B^5&=&x^5-5x^3\sqrt[n]{}\beta+5x\sqrt[n]{}\beta^2.
\end{eqnarray*}
and finally
\[
\begin{split}
&\sqrt[n]{}A^n+\sqrt[n]{}B^n\\
&=x^n-nx^{n-2}\sqrt[n]{}\beta+\frac{n(n-3)}{1\cdot 2}x^{n-4}\sqrt[n]{}\beta^2\\
&-\frac{n(n-4)(n-5)}{1\cdot 2\cdot 3}x^{n-6}\sqrt[n]{}\beta^3
+\frac{n(n-5)(n-6)(n-7)}{1\cdot 2\cdot 3\cdot 4}x^{n-8}\sqrt[n]{}\beta^4
-\textrm{etc.}\\
&=\alpha.
\end{split}
\]
This is the equation to be resolved whose resolvent is
\[
z^{n-1}=\alpha z^{n-2}-\beta z^{n-3} \quad \textrm{or} \quad z^2=\alpha z-\beta.
\]

\S 10. Moreover, not only in this way
will one root of the equation
\[
x^n-nx^{n-2}\sqrt[n]{}\beta+\frac{n(n-3)}{1\cdot 2}x^{n-4}\sqrt[n]{}\beta^2
-\textrm{etc.}=\alpha
\]
be found,
\[
x=\sqrt[n]{}A+\sqrt[n]{}B,
\]
but any other will also satisfy it
\[
x=\mu\sqrt[n]{}A+\nu\sqrt[n]{}B,
\]
providing $\mu^n=\nu^n=\mu\nu=1$, 
because it can be made in $n$ different ways.
So if $n=5$, there will be five roots of the equation
\[
x^5-5x^3\sqrt[5]{}\beta+5x\sqrt[5]{}\beta^2=\alpha
\]
as follow:

I. $x=\sqrt[5]{}A+\sqrt[5]{}B$,

II. $x=\frac{-1-\surd 5+\surd(-10+2\surd 5)}{4}\sqrt[5]{}A
+\frac{-1-\surd 5-\surd(-10+2\surd 5)}{4}\sqrt[5]{}B$,

III. $x=\frac{-1-\surd 5-\surd(-10+2\surd 5)}{4}\sqrt[5]{}A+
\frac{-1-\surd 5+\surd(-10+2\surd 5)}{4}\sqrt[5]{}B$,

IV. $x=\frac{-1+\surd 5+\surd(-10-2\surd 5)}{4}\sqrt[5]{}A
+\frac{-1+\surd 5-\surd(-10-2\surd 5)}{4}\sqrt[5]{}B$,

V. $x=\frac{-1+\surd 5-\surd(-10-2\surd 5)}{4}\sqrt[5]{}A+
\frac{-1+\surd 5+\surd(-10-2\surd 5)}{4}\sqrt[5]{}B$.

For all these coefficients are surdsolid roots of unity, and
the product by joining two is $=1$.

In a similar way beside unity itself, there are six roots
of the seventh power
of unity,
of which three pairs produce unity by multiplication, which are the six
roots of this equation
\[
y^6+y^5+y^4+y^3+y^2+y^1+1=0.
\]
However for finding these, 
the whole task is the resolution of a cubic equation;
for in fact, each equation of the sixth degree of the form
\[
y^6+ay^5+by^4+cy^3+by^2+ay+1=0,
\]
which does not change by putting $\frac{1}{y}$ in place of $y$,
can be resolved by means of cubic equations. 
Since since is often useful for finding roots, I will explain it briefly.

\S 11. These kind of equations, which when
$\frac{1}{y}$ is put in place of $y$ do not change their form,
I call {\em reciprocal}. These, if the maximum
dimension of $y$ is an odd number, can always be divided by $y+1$ and 
the resulting equation will also be reciprocal, in which the maximum
dimension of $y$ will be even.
It therefore suffices to have only considered equations of even dimensions, and to
have revealed the method for resolving them.

Thus first let this equation of the fourth dimension be given
\[
y^4+ay^3+by^2+ay+1=0;
\]
let this be produced from two quadratics
\[
y^2+\alpha y+1=0
\]
and
\[
y^2+\beta y+1=0.
\]
With this multiplied out it becomes
\[
\alpha+\beta=a \quad \textrm{and} \quad \alpha\beta+2=b \quad
\textrm{or} \quad \alpha\beta=b-2.
\]
Therefore $\alpha$ and $\beta$ will be the two roots of this equation
\[
u^2-au+b-2=0
\]
and by this rule the four roots of the given equation will be revealed
by means of just quadratic equations.

Let the reciprocal equation of the sixth power
\[
y^6+ay^5+by^4+cy^3+by^2+ay+1=0
\]
be produced from these three quadratics
\[
y^2+\alpha y+1=0,
\]
\[
y^2+\beta y+1=0
\]
and
\[
y^2+\gamma y+1=0.
\]
Then it will become
\[
\alpha+\beta+\gamma=a,
\]
\[
\alpha \beta+\alpha \gamma+\beta \gamma=b-3
\]
and
\[
\alpha\beta\gamma=c-2\alpha-2\beta-2\gamma=c-2a.
\]
Therefore $\alpha,\beta$ and $\gamma$ will be the three roots of this
cubic equation
\[
u^3-au^2+(b-3)u-c+2a=0.
\]

Similarly, the reciprocal equation of the eighth power
\[
y^8+ay^7+by^6+cy^5+dy^4+cy^3+by^2+ay+1=0
\]
is produced from four quadratic equations
\[
y^2+\alpha y+1=0,
\]
\[
y^2+\beta y+1=0,
\]
\[
y^2+\gamma y+1=0
\]
and
\[
y^2+\delta y+1=0,
\]
which which it follows
\[
\alpha+\beta+\gamma+\delta=a,
\]
\[
\alpha\beta+\alpha\gamma+\alpha\delta+\beta\gamma+\beta\delta+\gamma\delta=b-4,
\]
\[
\alpha\beta\gamma+\alpha\beta\delta+\alpha\gamma\delta+\beta\gamma\delta=c-3a
\]
and
\[
\alpha\beta\gamma\delta=d-2b+2.
\]
Therefore the coefficients $\alpha,\beta,\gamma,\delta$ are the four roots
of this equation
\[
u^4-au^3+(b-4)u^2-(c-3a)u+d-2b+2=0.
\]

The equation of the tenth order
\[
y^{10}+ay^9+by^8+cy^7+dy^6+ey^5+dy^4+cy^3+by^2+ay+1=0
\]
will be produced from five of these
\begin{eqnarray*}
y^2+\alpha y+1&=&0,\\
y^2+\beta y+1&=&0,\\
y^2+\gamma y+1&=&0,\\
y^2+\delta y+1&=&0,\\
y^2+\epsilon y+1&=&0,
\end{eqnarray*}
in which $\alpha,\beta,\gamma,\delta,\epsilon$ are the five roots of this
equation
\[
u^5-au^4+(b-5)u^3-(c-4a)u^2+(d-3b+5)u-e+2c-2a=0.
\]

And in general the reciprocal
equation
\[
\begin{split}
&y^{2n}+ay^{2n-1}+by^{2n-2}+cy^{2n-3}+dy^{2n-4}+ey^{2n-5}+fy^{2n-6}+
\cdots+py^n+\cdots\\
&+fy^6+ey^5+dy^4+cy^3+by^2+ay+1=0
\end{split}
\]
will be resolved into $n$ quadratic equations
\begin{eqnarray*}
y^2+\alpha y+1&=&0,\\
y^2+\beta y+1&=&0,\\
y^2+\gamma y+1&=&0,\\
y^2+\delta y+1&=&0\\
\textrm{etc.}&&
\end{eqnarray*}
And the coefficients $\alpha,\beta,\gamma,\delta$ etc. will be
the roots of this equation of $n$ dimensions
{\tiny
\[
\begin{array}{llllll}
u^n&-au^{n-1}&+bu^{n-2}&-cu^{n-3}&+du^{n-4}&-eu^{n-5}\\
&&-n&+(n-1)a&-(n-2)b&+(n-3)c\\
&&&&+\frac{n(n-3)}{1\cdot 2}&-\frac{(n-1)(n-4)}{1\cdot 2}a\\
+fu^{n-6}&-gu^{n-7}&+hu^{n-8}&-\textrm{etc.}&&\\
-(n-4)d&+(n-5)e&-(n-6)f&&&\\
+\frac{(n-2)(n-5)}{1\cdot 2}b&-\frac{(n-3)(n-6)}{1\cdot 2}c&+\frac{(n-4)(n-7)}{1\cdot 2}d&&&\\
-\frac{n(n-4)(n-5)}{1\cdot 2\cdot 3}&+\frac{(n-1)(n-5)(n-6)}{1\cdot 2\cdot 3}a&-
\frac{(n-2)(n-6)(n-7)}{1\cdot 2\cdot 3}b\\
&&+\frac{n(n-5)(n-6)(n-7)}{1\cdot 2\cdot 3\cdot 4}&&&\\
=0.&&&&&
\end{array}
\]
}

\S 12. Because the end term of each quadratic equation dividing the given equation is unity, it is clear that the product of the pairs of roots of the given
equation are unity.
 Thus joining together these pairs of terms
$\sqrt[n]{}A$ and $\sqrt[n]{}B$,
all the all the roots of the proposed equation
\S 9 may be obtained. 

\S 13. If in a reciprocal equation all the terms beside the extremes
and the middle are absent, as in
\[
y^{2n}+py^n+1=0,
\]
its divisors
\[
\begin{split}
&y^2+\alpha y+1,\\
&y^2+\beta y+1,\\
&y^2+\gamma y+1\\
&\textrm{etc.}
\end{split}
\]
will be obtained by substituting for $\alpha,\beta,\gamma,\delta$ etc.
the roots of this equation
\[
u^n-nu^{n-2}+\frac{n(n-3)}{1\cdot 2}u^{n-4}-\frac{n(n-4)(n-5)}{1\cdot 2\cdot 3}u^{n-6}
+\cdots \pm p=0,
\]
where $+p$ should be taken if $n$ is an even number, at $-p$ if
$n$ is odd.\footnote{Translator: Rudio notes that in fact if $n$ is even
then the last term should be $p \pm 2$.} From this it is apparent that
this equation agrees with the equation 
\[
x^n-nx^{n-2}\sqrt[n]{}\beta+\cdots=\alpha
\]
resolved in
\S 9, and from this all the divisors can be assigned.

\S 14. The above mentioned resolution of the formula $y^{2n}+py^n+1$ into factors
has a great use in the integration of the differential
formula\footnote{Translator: That is, we factor the denominator into
quadratic factors, and then use
partial fractions.}
\[
\frac{dy}{y^{2n}+py^n+1}
\]
which has already been much treated by Geometers. For
with the denominator resolved into its factors $y^2+\alpha y+1,y^2+\beta y+1$ etc., the entire integration is reduced to the quadrature of the circle
or the hyperbola.
It is further very helpful that the equation
\[
u^n-nu^{n-2}+\frac{n(n-3)}{1\cdot 2}u^{n-4}-\cdots \pm p=0,\footnote{Translator: See previous note.}
\]
from which $\alpha,\beta,\gamma$ etc. are determined, 
involves the division of the arcs of a circle into $n$ parts, and thus the
coefficients $\alpha,\beta,\gamma$ etc. may be easily found.

\S 15. However let us return to the question of eliciting the equations
to be resolved from the resolvent equations.
And let the resolvent equation be
\[
z^3=\alpha z^2-\beta z+\gamma,
\]
whose three roots are $A,B,C$; it will therefore be
\[
\alpha=A+B+C, \quad \beta=AB+AC+BC \quad \textrm{and} \quad \gamma=ABC.
\]
And thus the root $x$ of the equation to be resolved will be
\[
=\sqrt[n]{}A+\sqrt[n]{}B+\sqrt[n]{}C.
\]
Let us put
\[
p=\sqrt[n]{}AB+\sqrt[n]{}AC+\sqrt[n]{}BC,
\]
and with this done it will be
\[
\sqrt[n]{}A^2+\sqrt[n]{}B^2+\sqrt[n]{}C^2=x^2-2p
\]
and
\[
\sqrt[n]{}A^2B^2+\sqrt[n]{}A^2C^2+\sqrt[n]{}B^2C^2=p^2-2x\sqrt[n]{}\gamma;
\]
and so on, so that it follows:
\[
\begin{split}
&\sqrt[n]{}A^3+\sqrt[n]{}B^3+\sqrt[n]{}C^3=x^3-3px+3\sqrt[n]{}\gamma,\\
&\sqrt[n]{}A^3B^3+\sqrt[n]{}A^3C^3+\sqrt[n]{}B^3C^3=p^3-3px\sqrt[n]{}\gamma+
3\sqrt[n]{}\gamma^2;
\end{split}
\]

\[
\begin{split}
&\sqrt[n]{}A^4+\sqrt[n]{}B^4+\sqrt[n]{}C^4=x^4-4px^2+4x\sqrt[n]{}\gamma+2p^2,\\
&\sqrt[n]{}A^4B^4+\sqrt[n]{}A^4C^4+\sqrt[n]{}B^4C^4=p^4-4p^2x\sqrt[n]{}\gamma+
4p\sqrt[n]{}\gamma^2+2x^2\sqrt[n]{}\gamma^2;
\end{split}
\]

\[
\begin{split}
&\sqrt[n]{}A^5+\sqrt[n]{}B^5+\sqrt[n]{}C^5=x^5-5px^3+5x^2\sqrt[n]{}\gamma+
5p^2x-5p\sqrt[n]{}\gamma,\\
&\sqrt[n]{}A^5B^5+\sqrt[n]{}A^5C^5+\sqrt[n]{}B^5C^5=p^5-5p^3x\sqrt[n]{}\gamma
+5p^2\sqrt[n]{}\gamma^2+5px^2\sqrt[n]{}\gamma^2-5x\sqrt[n]{}\gamma^3.
\end{split}
\]

It is easily seen how this table can be continued further. Namely it is
\[
\begin{split}
&\sqrt[n]{}A^m+\sqrt[n]{}B^m+\sqrt[n]{}C^m=x(\sqrt[n]{}A^{m-1}+
\sqrt[n]{}B^{m-1}+\sqrt[n]{}C^{m-1}\\
&-p(\sqrt[n]{}A^{m-2}+\sqrt[n]{}B^{m-2}+\sqrt[n]{}C^{m-2})+
\sqrt[n]{}\gamma(\sqrt[n]{}A^{m-3}+\sqrt[n]{}B^{m-3}+
\sqrt[n]{}C^{m-3})
\end{split}
\]
and
\[
\begin{split}
&\sqrt[n]{}A^mB^m+\sqrt[n]{}A^mC^m+\sqrt[n]{}B^mC^m
=p(\sqrt[n]{}A^{m-1}B^{m-1}+\sqrt[n]{}A^{m-1}C^{m-1}+
\sqrt[n]{}B^{m-1}C^{m-1}\\
&-x\sqrt[n]{}\gamma(\sqrt[n]{}A^{m-2}B^{m-2}+
\sqrt[n]{}A^{m-2}C^{m-2}+\sqrt[n]{}B^{m-2}C^{m-2})\\
&+\sqrt[n]{}\gamma^2(\sqrt[n]{}A^{m-3}B^{m-3}+
\sqrt[n]{}A^{m-3}C^{m-3}+\sqrt[n]{}B^{m-3}C^{m-3}).
\end{split}
\]

\S 16. I have also observed other not negligible properties of these progressions.
For by putting
\[
\sqrt[n]{}A^m+\sqrt[n]{}B^m+\sqrt[n]{}C^m=R
\]
and
\[
\sqrt[n]{}A^mB^m+\sqrt[n]{}A^mC^m+\sqrt[n]{}B^mC^m=S
\]
it will be
\[
\sqrt[n]{}A^{2m}+\sqrt[n]{}B^{2m}+\sqrt[n]{}C^{2m}=R^2-2S
\]
and
\[
\sqrt[n]{}A^{2m}B^{2m}+\sqrt[n]{}A^{2m}C^{2m}+\sqrt[n]{}B^{2m}C^{2m}=
S^2-2R\sqrt[n]{}\gamma^m.
\]
In a similar way it is also
\[
\sqrt[n]{}A^{3m}+\sqrt[n]{}B^{3m}+\sqrt[n]{}C^{3m}=R^3-3RS+3\sqrt[n]{}\gamma^m
\]
and
\[
\sqrt[n]{}A^{3m}B^{3m}+\sqrt[n]{}A^{3m}C^{3m}+\sqrt[n]{}B^{3m}C^{3m}=S^3-
3RS\sqrt[n]{}\gamma^m+3\sqrt[n]{}\gamma^{2m}.
\]
And these series proceed in this way as straightforwardly as
the preceding.

\S 17. If $n=2$, it will be
\[
\alpha=x^2-2p \quad \textrm{and} \quad \beta=p^2-2x\surd \gamma
\]
and by joining these two equations one will have
\[
x=\surd A+\surd B+\surd C
\]
and
\[
p=\surd AB+\surd AC+\surd BC;
\]
also, $A,B$ and $C$ are the three roots of this cubic equation
\[
z^3=\alpha z^2-\beta z+\gamma.
\]
Thus eliminating the letter $p$ from these two equations yields
\[
\Big(\frac{x^2-\alpha}{2}\Big)^2-2x\surd y=\beta
\]
or
\[
x^4-2\alpha x^2-8x\surd \gamma=4\beta-\alpha^2,
\]
and a root $x$ of this equation is known, of course $=\surd A+\surd B+\surd C$;
this equation is consistent with that which was resolved in \S 5.

In a similar way, when two equations of the kind
\[
x^3-3px+3\sqrt[3]{}\gamma=\alpha
\]
and
\[
p^3-3px\sqrt[3]{}\gamma+3\sqrt[3]{}\gamma^2=\beta,
\]
occur, it will be
\[
x=\sqrt[3]{}A+\sqrt[3]{}B+\sqrt[3]{}C
\]
and
\[
p=\sqrt[3]{}AB+\sqrt[3]{}AC+\sqrt[3]{}BC
\]
with $A,B$ and $C$ being roots of the equation
\[
z^3=\alpha z^2-\beta z+\gamma
\]
as before. Or with the letter $p$ eliminated, an equation between $x,\alpha,\beta,\gamma$ follows, whose root $x$ will become known.

In exactly the same way, with these two equations 
\[
x^4-4px^2+4x\sqrt[4]{}\gamma+2p^2=\alpha
\]
and
\[
p^4-4p^2x\sqrt[4]{}\gamma+4p\surd \gamma+2x^2\surd \gamma=\beta
\]
occurring, it will be
\[
x=\sqrt[4]{}A+\sqrt[4]{}B+\sqrt[4]{}C
\]
and
\[
p=\sqrt[4]{}AB+\sqrt[4]{}AC+\sqrt[4]{}BC
\]
and again $A,B$ and $C$ are the roots of this equation
\[
z^3=\alpha z^2-\beta z+\gamma.
\]
To easily eliminate $p$, let us put
\[
x^2-2p=R \quad \textrm{and} \quad p^2-2x\sqrt[4]{}\gamma=S
\]
and it will be
\[
R^2-2S=\alpha \quad \textrm{and} \quad S^2-2R\surd\gamma=\beta.
\]
Now having eliminated $p$ from these two equations, we will have
\[
x^4=2Rx^2+8x\sqrt[4]{}\gamma+4S-R^2.
\]
Let us compare this equation with the original
\[
x^4=ax^2+bx+c;
\]
it will be
\[
R=\frac{a}{2}, \sqrt[4]{}\gamma=\frac{b}{8} \quad \textrm{or} \quad
\gamma=\frac{b^4}{4096} \quad \textrm{and} \quad S=\frac{c}{4}+\frac{a^2}{16}.
\]
Hence we will thus obtain
\[
\alpha=\frac{a^2}{8}-\frac{c}{2} \quad \textrm{and} \quad \beta=\frac{c^2}{16}+
\frac{a^2c}{32}+\frac{a^4}{256}-\frac{ab^2}{64}.
\]
Whence $A,B$ and $C$ will be the three roots of this equation
\[
z^3=\Big(\frac{a^2}{8}-\frac{c}{2}\Big)z^2-\Big(\frac{a^4}{256}+\frac{a^2c}{32}+\frac{c^2}{16}
-\frac{ab^2}{64}\Big)z+\frac{b^4}{4096},
\]
which agrees remarkably with that found in
\S 7.

\S 18. Therefore whenever it happens that the calculation leads to two
equations involving two unknowns $x$ and $p$, as appear in the formulas
in \S 15,
 the value of either can be assigned even if having eliminated the other
the equation becomes very complicated. Thus in these cases it will be expedient
to do the calculation not with one equation and one unknown, but to retain
two equations involving two unknowns
and investigate whether perhaps they might be contained among those
formulas,
which I have become convinced can often happen if the calculation
is done correctly.

\S 19. Like how we treated the resolvent equations
\[
z^2=\alpha z-\beta \quad \textrm{and} \quad z^3=\alpha z^2-\beta z+\gamma,
\]
the equation
\[
z^4=\alpha z^3-\beta z^2+\gamma z-\delta
\]
should be moved to next and treated in the same way.
Namely if its roots are $A,B,C$ and $D$,
put
\[
\sqrt[n]{}A+\sqrt[n]{}B+\sqrt[n]{}C+\sqrt[n]{}D=x
\]
and
\[
\sqrt[n]{}AB+\sqrt[n]{}AC+\sqrt[n]{}AD+\sqrt[n]{}BC+\sqrt[n]{}BD
+\sqrt[n]{}CD=p
\]
and also
\[
\sqrt[n]{}ABC+\sqrt[n]{}ABD+\sqrt[n]{}ACD+\sqrt[n]{}BCD=q
\]
and then look for expressions for
\[
\sqrt[n]{}A^m+\sqrt[n]{}B^m+\textrm{etc.} \quad \textrm{and}
\quad \sqrt[n]{}A^mB^m+\sqrt[n]{}A^mC^m+\textrm{etc.}
\quad \textrm{and for} \quad \sqrt[n]{}A^mB^mC^m+\textrm{etc.}
\]
With this done, three equations containing $x,p$ and $q$ will always be found
for any value of $m$. And in a similar way the three unknowns occurring in these
three equations will be determined.

\S 20. I suspect that by putting
\[
x=\sqrt[5]{}A+\sqrt[5]{}B+\sqrt[5]{}C+\sqrt[5]{}D,
\]
a rational equation can be constructed in which $x$ has no more than $5$ dimensions,
even if this seems usually impossible. Namely, as in \S 17 from the equations
\[
x^4-4px^2+4x\sqrt[4]{}\gamma+2p^2=\alpha
\]
and
\[
p^4-4p^2x\sqrt[4]{}\gamma+4p\surd\gamma+2x^2\surd\gamma=\beta
\]
in eliminating $p$ we obtained an equation with not more than $4$ dimensions,
which would seem equally hard to happen, so perhaps by a similar use for
the fifth power an artifice could be come to such that the equation
\[
x^5=ax^3+bx^2+cx+d
\]
could finally be resolved.
In my opinion, what is most important in completing this is that
$\alpha,\beta,\gamma,\delta$ be determined from $a,b,c,d$ and not the other way
around;
for in that case the equation would rise to a much higher power
than is useful.
However I leave the completion of this problem to others who delight in
these occupations, or to myself at another time,
and I am contented to have perhaps shown a suitable and natural approach.

\end{document}